\documentclass[12pt,a4paper]{amsart}
\usepackage[colorlinks=true,          
            linkcolor=blue,
            citecolor=red,
            urlcolor=blue]{hyperref}
\usepackage{a4wide,amssymb,amsfonts,amsmath,hyperref,mathrsfs,latexsym,stmaryrd,graphicx,mathrsfs,amsaddr}
\DeclareMathAlphabet{\mathpzc}{OT1}{pzc}{m}{it}

 \newtheorem{thm}{Theorem}[section]

\input xy
\xyoption{all}

\newcommand {\CC}{\mathbb{C}}

\newcommand {\PP}{\mathbb{P}}

\newcommand {\RR}{\mathbb{R}}

\newcommand {\HH}{\mathbb{H}}

\newcommand{\ZZ}{\mathbb{Z}}

\newcommand {\br}{{\bf r}}

\newcommand {\bB}{{\bf B}}
\newcommand {\bD}{{\bf D}}

\newcommand {\bF}{{\bf F}}
\newcommand {\bG}{{\bf G}}

\newcommand {\bO}{{\bf O}}

\newcommand {\bU}{{\bf U}}
\newcommand {\bV}{{\bf V}}

\newcommand{\rk}{\textrm{rk }}

\newcommand{\cB}{\mathcal{B}}
\newcommand{\cC}{\mathcal{C}}

\newcommand{\cF}{\mathcal{F}}
\newcommand{\cG}{\mathcal{G}}
\newcommand{\cH}{\mathcal{H}}

\newcommand{\cL}{\mathcal{L}}

\newcommand {\cO}{\mathcal{O}}
\newcommand{\cP}{\mathcal{P}}

\newcommand{\cX}{\mathcal{X}}

\newcommand{\cR}{\mathcal{R}}
\newcommand{\cS}{\mathcal{S}}
\newcommand{\cT}{\mathcal{T}}
\newcommand{\cU}{\mathcal{U}}
\newcommand{\cV}{\mathcal{V}}

\newcommand{\scB}{\mathscr{B}}

\newcommand{\scL}{\mathscr{L}}

\newcommand{\scU}{\mathscr{U}}

\newcommand{\scZ}{\scZ}

\newcommand {\fU}{\mathfrak{U}}

\newcommand{\fg}{\mathfrak{g}}

\newcommand{\ft}{\mathfrak{t}}

\newcommand{\fD}{\mathfrak{D}}

\newcommand{\mhom}{\textrm{Hom}}

\newcommand {\ad}{\textrm{ad } }

\newcommand{\spec}{\textrm{Spec }}

\newcommand{\tot}{\textrm{tot }}

\newcommand{\ctimes}{\otimes_\CC}
\newcommand{\ztimes}{\otimes_\ZZ}

\newcommand {\io}{\iota}
\newcommand {\fii}{\varphi}
\newcommand{\hookr}{\hookrightarrow}

\newcommand{\ram}{\textrm{Ram }}
\newcommand {\bra}{\textrm{Bra }}
\newcommand{\ev}{\textrm{ev}}
\newcommand{\sym}{\textrm{Sym}}
\newcommand{\kd}{L}
\newcommand{\diag}{\textrm{diag}}
\newcommand{\gr}{\textrm{gr}}

\newcommand{\Higgs}{{\bf Higgs}}

\newcommand{\Prym}{{\bf Prym}}

\newcommand{\Pic}{{\bf Pic}}
\newcommand{\rts}{\sf{root}}
\newcommand{\wts}{\sf{weight}}
\newcommand{\crts}{\sf{coroot}}
\newcommand{\cwts}{\sf{coweight}}
\newcommand{\chr}{\sf{char}}
\newcommand{\cchr}{\sf{cochar}}

\include{biblio}
\title[Donagi--Markman cubic for the generalised Hitchin system]{Donagi--Markman cubic for the generalised \\[5pt] Hitchin system}

\author{Ugo Bruzzo$^{\P\S}$ and Peter Dalakov$^{\P\ddag}$}

\address{\small $^\P$Scuola Internazionale Superiore di Studi Avanzati (SISSA),\\
Via Bonomea 265, 34136 Trieste, Italia\\[3pt]     
$^\S$Istituto Nazionale di Fisica Nucleare, Sezione di Trieste \\[3pt] 
$^\ddag$On leave from IMI-BAS, Sofia, Bulgaria}

\date{\today} 
\subjclass[2010]{14D20, 14D07,  14H70 }
\thanks{U.B.'s research is partly supported by PRIN ``Geometria delle variet\`a algebriche"  and by GNSAGA-INDAM. He  is a member of the VBAC  group. P.D.~is supported 
by a CERES--CEI research fellowship 
FP-7, Marie Curie Actions.}
\newtheorem*{thma}{Theorem A}
\newtheorem*{thmb}{Theorem B}

\newtheorem{proposition}[thm]{Proposition}

\theoremstyle{remark}

\theoremstyle{definition}

\begin{document}

\begin{abstract}
Donagi and Markman (1993) have shown that the infinitesimal period map for an algebraic completely integrable Hamiltonian system (ACIHS) is encoded
in a section of the third symmetric power of the cotangent bundle to the base of the system.
 For the ordinary Hitchin system the cubic
is given by a formula of Balduzzi and Pantev. We show that
the Balduzzi--Pantev formula  holds
on maximal rank symplectic leaves  of the $G$-generalised
Hitchin system.
\end{abstract}
\maketitle 
\setcounter{tocdepth}{1}
\tableofcontents 
\thispagestyle{empty}
\section{Introduction}
Algebraic completely integrable Hamiltonian systems (ACIHS) form a very special class of quasi-projective algebraic varieties. They
carry both a holomorphic symplectic (or Poisson) structure and the structure of a Lagrangian fibration, whose generic fibres are 
principal homogeneous spaces for
 abelian varieties (algebraic tori).  The varying Hodge structure on the Lagrangian tori allows one to write, after some choices, a holomorphic period map from the 
base of the fibration to the Siegel upper-half space.
Such structures often arise after complexification of real integrable systems appearing in
classical physics.

In the $C^\infty$-category any smooth family of tori admits locally (on the base) the structure of a Lagrangian fibration,
provided the dimension of the  base is  half the dimension of the total space. In the holomorphic category, 
 as discovered in \cite{donagi_markman_cubic}, for a Lagrangian structure to exist, the differential of the period map must be a
section of the third symmetric power of the cotangent bundle to the base, \emph{the Donagi--Markman cubic}.
Locally on the base,
the cubic 
 is given (in special coordinates) by the third derivatives of a holomorphic function $\bF$, called \emph{a holomorphic pre-potential}, while
the (marked) period map is  given by its Hessian $\textrm{Hess}(\bF)$. For example, if the ACIHS is the family of intermediate Jacobians over the moduli space of
gauged Calabi--Yau threefolds, the cubic is the Yukawa cubic (\cite{donagi_markman_cubic}).

In \cite{hitchin_sd}, \cite{hitchin_sb} N.Hitchin introduced a beautiful ACIHS, supported on a partial compactification
of the cotangent bundle to the coarse moduli space of  holomorphic $G$-bundles on a  Riemann surface $X$
of genus at least two. Its points correspond to (semi-stable) Higgs bundles on $X$. These are  pairs
 $(E,\theta)$ of a  $G$-bundle $E$ and  a holomorphic section of
$\ad E\otimes K_X$, where $G$ is   a semi-simple or reductive
complex Lie group.
 The ``conserved momenta''   are the spectral invariants of $\theta$,
 and the abelian varieties
  are Jacobians or Pryms of certain covers of $X$. 

In \cite{markman_thesis}, \cite{bottacin} the authors considered  natural generalisations of the 
Hitchin system, supported on    moduli spaces of meromorphic, i.e., $K_X(D)$-valued, Higgs bundles,
for  a divisor $D\geq 0$.  Such moduli spaces exist in any genus, provided $D$ is sufficiently positive,
and 
Markman and Bottacin showed that they carry a natural 
holomorphic Poisson structure. 

Many known examples of ACIHS appear  as special cases of this \emph{generalised Hitchin system}, see
the survey \cite{donagi_markman}.
In genus zero such examples are the 
 finite-dimensional coadjoint orbits
in   loop algebras and
 Beauville's (\cite{beauville_matrices}) polynomial matrix systems,
including geodesic flows on ellipsoids  and Euler--Arnold systems.
 In genus one we have elliptic Gaudin, Sklyanin (\cite{hurtubise_markman_sklyanin})
and Calogero--Moser (\cite{hurtubise_markman_cal_mo}) systems.

For any ACIHS, $\textrm{Im }\textrm{Hess}(\bF)$,
being symmetric and positive-definite,  provides a K\"ahler
metric on the base. Together with the family of
lattices, this data can be packaged into the so-called \emph{(integral) special K\"ahler structure}
(\cite{danfreed}), used by physicists to describe  massive vacua in global $N=2$ 
super-symmetry in 4 dimensions, e.g.,  in Seiberg--Witten theory. More specifically, Donagi and Witten
(\cite{donagi_witten}, \cite{donagi_swis}, \cite{markman_sw}) have proposed certain special cases of the
generalised Hitchin system as candidates for the Coulomb branch of the moduli of
vacua for 4-dimensional $N=2$ Yang--Mills theory with structure group $G_c$ and massive vector multiplets.

A surprising relation between $\bB$-model large $N$ duality
and ADE  Hitchin systems has been discovered in \cite{geom_trans}, \cite{ddp}.
 There the authors construct  
a family of quasi-projective Calabi--Yau 3-folds, whose family of intermediate Jacobians is isogenous to the
family of Hitchin Pryms. In that setup the Donagi--Markman cubic plays the role of the Yukawa cubic.

 Identifying the cubic is thus a natural and interesting question.
 For the ordinary Hitchin system ($D=0$) and $G=SL_2$ such a  formula
appears in \cite{geom_trans}, (47).  
Building on
 unpublished notes of T.Pantev,
Balduzzi (\cite{balduzzi}) has derived a formula  for arbitrary semi-simple $G$ and $D=0$, and
another variant of the  proof is sketched in \cite{hertling_hoev_posthum}.
Finally,   \cite{hoevenaars} treats the case of
 the Neumann system, i.e.,  $X=\PP^1$, $G=SL_2$ and $D=2\cdot\infty + \sum_i n_i q_i$,
using special properties of plane curves. In all of these examples the cubic is given by
a ``logarithmic derivative'' of a discriminant.

Our goal in this note is to show that the Balduzzi--Pantev formula holds 
along the (good) symplectic leaves of the generalised Hitchin system.
     \begin{thma}
Let $X$ be a smooth curve, $G$ a semi-simple complex algebraic Lie group, and
 $D\geq 0$  a  divisor on $X$ with  $L^2=K_X(D)^2$  very ample.
Let $o\in\cB\simeq H^0(X,\bigoplus L^{d_i})$    be a generic point, corresponding to a smooth
cameral cover $\pi_o:\widetilde{X}_o\to X$ with simple   ramification.
 Denote by $\fD:\cB\times X\to \tot L^{|\cR|}$  the discriminant, and by $\bB$ the set
\[
 \bB = \left(\{o\}+ H^0(X,\bigoplus_i L^{d_i}(-D))\right)\cap \scB\subset \cB,
\]
where $\scB\subset \cB$ is the locus of generic cameral covers.
Then the Donagi--Markman cubic for the (maximal rank) symplectic leaf  $ \left. \cS\right|_{\bB}\to \bB$
of the $G$-generalised Hitchin system 
is  
\[
 c_o: H^0(\widetilde{X}_o, \ft\ctimes K_{\widetilde{X}_o})^W \longrightarrow \sym^2 \left(H^0(\widetilde{X}_o, \ft\ctimes K_{\widetilde{X}_o})^W \right)^\vee
\]
\[
 c_o(\xi)(\eta,\zeta)= \frac{1}{2}\sum_{p\in\ram(\pi_o)} \textrm{Res}_p^2\ \left( \pi_o^\ast \left.\frac{\cL_{Y_\xi}(\fD)}{\fD}\right|_{\{o\}\times X}\eta\cup\zeta\right).
\]
Here 
$Y_\xi$ is the preimage of $\xi$ under the isomorphism
$T_{\bB,o}\simeq H^0(\widetilde{X}_o,\ft\ctimes K_{\widetilde{X}_o})^W $ induced by the canonical meromorphic 2-form on $\tot \ft\ctimes L$, and
$\cL_{Y_\xi}$ denotes Lie derivative.
      \end{thma}
    \begin{thmb}
With the same notation as above,
\[
 c_o(\xi,\eta,\zeta) =\sum_{p\in\ram(\pi_o)} \textrm{Res}^2_p\sum_{\alpha\in\cR}\frac{\alpha(\xi)\alpha(\eta)\alpha(\zeta)}{\alpha(\lambda_o)},
\]
where $\lambda_o\in H^0(\widetilde{X}_o,\ft\ctimes L)$ is the tautological section.
    \end{thmb}
The proof is   a local calculation for the universal family of cameral covers.
The main ingredient, making this possible, is  \emph{abelianisation} (\cite{don-gaits}),
which has a long history,  see
\cite{hitchin_sb,bnr,faltings,scognamillo_elem,donagi_spectral_covers,don-gaits}. The family of Albanese varieties, associated to the symplectic leaves of the Hitchin
system is a family of generalised Prym varieties for a family of cameral curves. By a theorem of Griffiths,
the differential of the period map of the latter is given by cup product with the Kodaira--Spencer class.
By  results of \cite{hurtubise_markman_rk2} and \cite{kjiri}, the symplectic structure
on the   cameral Pryms agrees with the one on the leaves of the Hitchin system, so the Kodaira--Spencer
class is the only missing ingredient.

The content of the paper is as follows. In Section \ref{hitchin} we review the generalised Hitchin
system and Markman's construction of the Poisson structure via symplectic reduction.
In Section \ref{VHS} we review period maps and the cubic condition. In Sections \ref{covers}
and \ref{abel}
we recall basics about cameral covers and  abelianisation.
We prove Theorems A and B  in Section \ref{proof}.
	  \subsection*{Acknowledgments} P.D.~thanks Tony Pantev and Calder Daenzer for helpful discussions related to this project,
and  SISSA  for  hospitality. %
	  \subsection*{Notation and conventions} 
We fix    a smooth, compact, connected Riemann surface  $X$ of genus $g$,  $D=\sum_{i=1}^{s} n_i q_i\geq 0$ 
a divisor. Set $L=K_X(D)$ and assume   $L^2$ is very ample.
   We fix $G=\bG(\CC)\supset B\supset T$: a semi-simple complex algebraic Lie group,  a
 Borel subgroup, and a maximal torus.
Then $W=N_G(T)/T$ is the  Weyl group,  $\cR\subset \ft^\vee$  the root system, and 
$\cR^+$ the positive roots.
 The 
 lattices $\rts_\fg\subset \chr_G\subset \wts_\fg$ and
$\cchr_\fg\subset \cwts_G\subset \cwts_\fg$ correspond to their names.
By  $\bG\br=\bG(n,V)$ we denote Grassmannian of $n$-dimensional subspaces of a vector space $V$ and by
 $Z(s)$ the zero locus of a section $s$.  The Hitchin base is $\cB=H^0(X,\ft\ctimes L/W)\simeq H^0(X,\oplus_i L^{d_i})$,
and $\scB=\cB\backslash \Delta$ is the complement of the discriminant locus. 
We denote  ramification and branch loci by $\ram$ and $\bra$, respectively.
%

 \medskip
      \section{The generalised Hitchin system}\label{hitchin}
	\subsection{Meromorphic Higgs bundles} 

A \emph{holomorphic  $\kd$-valued $G$-Higgs bundle} on $X$ is  a pair $(E,\theta)$, consisting of
 a (holomorphic) principal $G$-bundle $E\to X$ and   $\theta\in H^0(X, \ad E\otimes \kd)$.
The pair is \emph{(semi-) stable}, if for any \emph{$\theta$-invariant} reduction  $\sigma: X\to E/P$ of
$E$ to a maximal parabolic subgroup $P\subset G$, $\deg \sigma^\ast T_{G/P}>0$ (resp. $\geq 0$). Here $T_{G/P}$ is the vertical bundle of
$E/P\to X$. This is  Ramanathan's
(semi-)stability \cite{ramanathan}, but tested on $\theta$-invariant reductions only.

By the argument in Theorem 4.3 of \cite{ramanathan}, there exists an analytic 
coarse moduli space $\Higgs_{G,D} =\coprod_{c\in \pi_1(G)}\Higgs_{G,D,c}$ 
of semi-stable 
pairs $(E,\theta)$, whose
connected components
are labelled by the topological
type  of  $E$. It is expected to be algebraic and separated  
by suitable modifications of the  constructions in \cite{moduli2} and \cite{nitsure}, 
see Remark 2.5 in \cite{hurtubise_markman_sklyanin}.
For  \emph{everywhere regular} $\theta$  one can 
   construct the 
moduli space using the spectral correspondence, see Section 5.4 of \cite{donagi_spectral_covers}. 
	\subsection{Poisson structure}  
If  nonempty,    $\Higgs_{G,D,c}$ are   Poisson  (\cite{markman_thesis}, \cite{bottacin}), and in fact
 are ACIHS via 
the \emph{Hitchin map} $h_c$. This is
 a proper morphism $\Higgs_{G,D,c}\to \cB$, mapping $[(E,\theta)]$ to the spectral invariants of $\theta$.
   The  \emph{ Hitchin base}  $\cB$ (Section \ref{cameral})
 is 
 (non-canonically) isomorphic to $\bigoplus_i H^0(X,\kd^{d_i})$, where
  $d_i$ are the degrees of the basic $G$-invariant polynomials on $\fg$.
 We recall next some  definitions
 from \cite{donagi_markman}, and 
 the structure of the symplectic leaves    from  \cite{markman_thesis} and \cite{markman_sw}.
	     A Poisson structure on a smooth (analytic or quasi-projective algebraic) variety $M$ is 
 a  bivector  
$\Pi\in H^0(M,\Lambda^2 T_M)$, endowing $\cO_M$ with the structure of a sheaf of Lie algebras via the Poisson
bracket $\{f,g\}:= (df\wedge dg) (\Pi)$. The Lie subalgebra of \emph{Casimir functions} is 
$\ker \left(H^0(\cO_M)\to H^0(T_M) \right)$,  $f\mapsto X_f=\{f,\ \}$.

Such a $\Pi$ determines a sheaf morphism  $\Psi: T^\vee_M\to T_M$, and, consequently,
a  stratification of $M$ by  \emph{submanifolds} $M_k$ ($k$-even) on which
$\rk\Psi=k$. The strata are (local-analytically) foliated by $k$-dimensional leaves $S$  on which $\left.\Pi\right|_S$ is   symplectic.
These leaves are the integral leaves of  $\Psi(T^\vee_M)\subset T_M$ and the level sets of
the Casimirs.
	     A subvariety $L\subset (M,\Pi)$ is called \emph{Lagrangian}, if $L\subset \overline{S}$ for some symplectic leaf $S$, and
$L\cap S$ is a Lagrangian subvariety of $S$.
We say that $(M,\Pi)$ is an \emph{algebraically completely integrable Hamiltonian system (ACIHS)}, if there exist
a smooth variety $\cB$,   a  closed subvariety $\Delta\varsubsetneq \cB$, and
 a proper
flat morphism $h: M\to\cB$,
such that   $\left. h\right|_{\cB\backslash \Delta}$ has   Lagrangian fibres, which are   abelian torsors.
	    
Suppose     a  connected algebraic group $\cG$    acts on $(M,\Pi)$.  
The dual of $\mathbf{g}=\textrm{Lie }\cG$
carries the 
Kostant--Kirillov Poisson structure, whose symplectic leaves 
are the coadjoint orbits
$\mathbf{O}\subset\mathbf{g}^\vee$.
The action is said to be  \emph{Poisson}, if the Lie-algebra  homomorphism $\mathbf{g}\to H^0(T_M)$  factors through
Lie-algebra homomorphisms
 $\mathbf{g}\to H^0(\cO_M)\to H^0(T_M)$, the second arrow being $f\mapsto X_f$. The first arrow corresponds to
a \emph{moment map}
$\mu: M\to \mathbf{g}^\vee$,  which is Poisson. 
If $\mu$ is submersive with connected fibres, and if  $M/\cG$ exists, then it carries a canonical (Marsden--Weinstein) Poisson structure, whose
 symplectic leaf through $m\in M$ is 
\[
 \mu^{-1}(\mathbf{O}_{\mu(m)})/\cG\simeq \mu^{-1}(\mu(m))/\cG_{\mu(m)}.
\]

Let $E$ be a  $G$-bundle, 
and $\eta: \left. E\right|_{D}\simeq D\times G=\spec (H^0(\cO_D)\otimes H^0(\cO_G))$
 a framing.
Let $\widetilde{G}_D=\widetilde{\bG}_D(\CC)$ be  the group of maps from $D$ to $G$. We recall that
$\widetilde{\bG}_D$ is
defined in terms of its functor of points, by requiring that
$\mhom(\spec A, \widetilde{\bG}_D)= \mhom (\spec A\otimes H^0(\cO_G),\bG)$, for any  $\CC$-algebra $A$. It  is the product of the groups of
$(n_i-1)$-jets
of maps  $D\to G$, i.e.,  $\widetilde{G}_D= \prod_{i=1}^s G_{n_i-1}$, where
$\mhom(\spec A, \widetilde{\bG}_n)\simeq \mhom (\spec A\otimes \CC[t]/t^{n+1},\bG)$. 
For example, an embedding $G\subset GL_N$ would identify  $G_n$ with the set  of $\CC[t]/t^{n+1} $-valued matrices,
satisfying the equations of $G\mod t^{n+1}$.
The \emph{level group} is $G_D= \widetilde{G}_D/Z(G)$, with   $Z(G)$  embedded diagonally.
We identify
$\fg_D^\vee=\widetilde{\fg}_D^\vee\simeq \fg^\vee \otimes H^0(\left.K(D)\right|_D)$
by the duality pairing 
\[\xymatrix@1{H^0(\cO_D)\otimes H^0(\left. K(D)\right|_D)\ar[r]& H^0(\left. K(D)\right|_D)\ar[r]^-{Res}&H^1(K)\simeq \CC }  .\]
	  The moduli space $\cP^{st}=\cP_{G,D,c}$ of (stable) framed $G$-bundles of topological type 
$c$ with level-$D$ structure parametrises   isomorphism classes of stable pairs $(E,\eta)$,
 (\cite{markman_thesis}). Its cotangent bundle 
$T^\vee\cP^{st}$ parametrises  classes of triples $(E,\eta,\theta)$,
$\theta\in H^0(X,\ad E\otimes K_X(D))$. The natural action of $G_D$
  on $\cP^{st}$  lifts  to $T^\vee\cP^{st}$,
and the canonical moment map for a lifted action is  
  $\mu: \left[(E,\eta,\theta)\right]\mapsto Ad_\eta (\left.\theta\right|_D)$,
see \cite{markman_thesis} 6.12, \cite{markman_sw} (10).
Markman identified   $T^\vee\cP^{st}/G_D$ with a dense open  $\Higgs^{sm}_{G,D,c}$ and extended  the Poisson structure
 to the partial compactification.
	    
Let  $\cB_0:= H^0\left(\oplus_i L^{d_i}(-D)\right)\subset \cB$. 
By Proposition 8.8 of \cite{markman_thesis},
there is a canonical  isomorphism
$ \cB/\cB_0 \simeq \fg^\vee_D \sslash G_D \simeq \CC^{r\deg D}$,
\[
  \bigoplus_{k=1}^r \frac{H^0(L^{d_k})}{H^0(L^{d_k}(-D))}\simeq \bigoplus_{i=1}^s \fg^\vee\otimes H^0(\left. K(D_i)\right|_{D_i})\sslash G_{n_i-1}  \simeq \bigoplus_{i=1}^s \CC^{rn_i} =
\CC^{r\deg D}
\]
sending  $[f]\in H^0(L^{d_k})/H^0(L^{d_k}(-D))$  to the coefficients of
its $(n_i-1)$-jet in a holomorphic trivialisation at $q_i\in \textrm{supp}(D)$, $\forall i=1..s$:
\[
 \left(a_0 +\ldots +a_{n_i-1}z^{n_i-1}+\ldots +a_{n_id_k-1} z^{n_id_k-1}\right)\frac{dz^{\otimes d_k}}{z^{n_i d_k}}+\ldots \longmapsto \left(a_0,\ldots, a_{n_i-1}\right)\in\CC^{n_i}.
\]
The following diagram is a variant of \cite{markman_thesis}, Proposition 8.8:
\[
 \xymatrix{   & T^\vee\cP^{st}\ar[dr]^-{\mu}\ar[dl]& \\
	      \Higgs^{sm}_{G,D,c}\ar[d]_-{h}\ar[dr]^-{\overline{h}} & & \fg_D^\vee\ar[d]\\
		\cB\ar[r]& \cB/\cB_0\ar[r]^-{\simeq}&\fg^\vee_D \sslash G_D\\  
}.
\]
The symplectic leaves  $S=\mu^{-1}(\bO)/G_D$ are  labelled by
coadjoint orbits $\bO=\prod_i \bO_i\subset \fg^\vee_D$. The symplectic foliation  is a refinement of the
foliation by fibres of $\overline{h}$, every $\overline{h}$-fibre containing a unique leaf of maximal
($=\dim M_{\Higgs}-r\deg D$) rank.


\medskip
	  \section{VHS and the cubic condition}\label{VHS}
    \subsection{Variation of Hodge structures}
Let $h: \cH\to \cB$ be a holomorphic family of polarised compact connected K\"ahler manifolds over
a complex manifold $\cB$. 
The  varying Hodge filtration on the first de Rham cohomology of the fibres $h^{-1}(b)=:\cH_b$ is  
the prototypical  example of an integral, weight-1 variation of polarised Hodge structures (VHS),  encoded in the quadruple 
 $(\cF^\bullet,   \cF_\ZZ, \nabla^{GM}, S)$.
Here $\cF=R^1 h_\ast \CC\otimes\cO_B$,  a 
 holomorphic vector bundle with
 fibres
 $\cF_b = H^1(\cH_b,\CC)$.  It carries a flat holomorphic connection 
 $\nabla^{GM}: \cF\to \cF\otimes \Omega^1_B$, a  $\nabla^{GM}$-flat subbundle of lattices $\cF_\ZZ$, a holomorphic flag
$\cF^1\subset \cF^0=\cF$ and a polarisation  $S$.

The Gauss--Manin connection $\nabla^{GM}$ can be defined 
 topologically, using  Ehresmann's theorem:  $\cH$ is a  fibre bundle, so
 homotopy invariance of de Rham cohomology implies that
$R^1 h_\ast \CC$ is a locally constant sheaf, and  the connection is  then  expressed  by a Cartan--Lie
formula (\cite{voisin1}).
A holomorphic description, due to  Katz and Oda (\cite{katz-oda}), is  as follows.
The complex   $\Omega^\bullet_\cH$ can be equipped with
the Koszul--Leray filtration 
$ L^i\Omega^\bullet_{\cH}=h^\ast\Omega^i_\cB\wedge \Omega^\bullet_{\cH}[-i]$.
 Then in the 
 associated
spectral sequence   
 $(E_1^{\bullet,0},d_1)$ is identified with  $(\Omega_\cB^\bullet(\cF),\nabla^{GM} )$, via
$ h^{-1}\cO_\cB\simeq _{qis}\Omega^\bullet_{\cH/\cB}[-1]$.

The Hodge filtration 
 $\cF^1=h_\ast\Omega^1_{\cH/\cB}\subset \cF^0=\cF$
is induced by 
the 
 stupid filtration  $\Omega^{\geq 1}_{\cH/\cB}[-1]\subset \Omega^\bullet_{\cH/\cB}$. 
We denote by $\gr \nabla^{GM}$ 
the $\Omega^1_B$-valued homomorphism $\cF^1\to \Omega^1_B\otimes \cF^0/\cF^1$, which is  Simpson's original 
Higgs field on $\cF^1\oplus \cF^0/\cF^1$.

Integral cohomology provides a subbundle of lattices $\cF_\ZZ=R^1 h_\ast \ZZ$ and the  induced 
 real structure splits  the Hodge filtration, giving  $\cF\simeq_{C^\infty}\cF^1 \oplus\overline{\cF^1} $. 

 The fibrewise
K\"ahler structure   provides a $\nabla^{GM}$-flat, non-degenerate, sesqui-linear, \emph{skew-symmetric}
 pairing $\xymatrix@1{S: \cF^0\otimes\overline{\cF^0}\ar[r]&\ \cC^\infty_\cB}$,
$ S_b([\eta],[\psi])=i\int_{\cH_b} \eta\wedge \overline{\psi}\wedge \omega_b^{\dim\cH/\cB-1},$
which  induces a $C^\infty$-isomorphism  $\overline{\cF^1}\simeq_{C^\infty}\cF^{1\ \vee} $. 
The Hodge--Riemann bilinear relations state that $1)$
 $H^1(\cH_b,\cO)$ and $H^0(\cH_b,\Omega^1)$ are Lagrangian for $S_b(\bullet,\bar\bullet)$ and $2)$  $S_b>0$
on 
$ H^0(\cH_b,\Omega^1)$.
	\subsection{The period map}
 The VHS data      are  equivalent to a period map from $\cB$ to a classifying space
(period domain), which for weight one is Siegel's upper half space $\HH$.
Let $o\in \cB$ be a base point and $\check \bD\subset \bG\br=\bG\br (h^{1,0}, H^1(\cH_o,\CC))$
 the subvariety of \emph{Lagrangian} (for $S_o$)
 subspaces, which is a homogeneous space
 $\textrm{Sp}(H^1(\cH_o,\CC))/\textrm{Stab}\left[ H^0(\Omega^1_{\cH_o})\subset H^1(\CC)\right]$.
The second bilinear relation determines  an open  subvariety  
$\bD\subset \check \bD$, isomorphic to  $\textrm{Sp}(2h^{1,0},\RR)/$ $ U(h^{1,0})$.
Parallel transport then  gives rise to a (set-theoretic) period map $\Phi: \cB\to \bD/\Gamma$,
where $\Gamma$ is the monodromy group of the VHS.
\begin{thm}[\cite{griffiths_periods}]
The period map $\Phi$  is  locally liftable and holomorphic.
\end{thm}
Hence on a
contractible neighbourhood  $\cU\subset(\cB,o)$
 there is a holomorphic lift 
\[
 \xymatrix@1{\widetilde{\Phi}:& \cU\ar[r]& \ \bD \subset \bG\br (h^{0,1}, H^1(\cH_o,\CC))}
\]
\[
 b\longmapsto \left[H^0(\cH_b,\Omega^1)\subset H^1(\cH_b,\CC)\simeq H^1(\cH_o,\CC)\right].
\]
We can make  $\widetilde{\Phi}$ more explicit as follows. The Hodge decomposition on
$H^1(\cH_o,\CC)$ gives a standard coordinate chart  $ H^0(\Omega^1)^\vee\otimes H^1(\cO)\simeq_S H^0(\Omega^1)^{\vee\ \otimes 2}\subset \bG\br$.
Then $\check \bD\cap H^0(\Omega^1)^{\vee\ \otimes 2} =\sym^2\left(H^0(\Omega^1)^\vee \right) $, and $\bD$ consists of the elements
with $\textrm{Im }>0$.

The classical way of explicating  this map is by choosing a \emph{marking}.
A harmonic basis of $H^0(\cH_o,\Omega^1)$ identifies $\bD$ with Siegel's upper-half
space $\HH_{h^{1,0}}$. 
Since $\cU\subset \cB$ is contractible,  we have 
$H_1(\cU,\ZZ)\simeq H_1(\cH_o,\ZZ)$. Choosing a basis of $H_1$ mod torsion,
with normalised A-periods, gives a period matrix of the form
 $[\Delta_\delta, Z]$, and  $\widetilde{\Phi}(b)=Z(b)\in\HH_{h^{1,0}}$. Here $\Delta_\delta=\diag(\delta_1,\ldots)$,
where $\delta_i$ are the polarisation divisors.
	\subsection{A theorem of Griffiths}
We can associate to the family $h$  various  linear objects:
the vertical bundle $\cV=h_\ast T_{\cH/\cB}$,   the relative 
Jacobian ($\overline{\cF^1}/\cF_{\ZZ}$) and
Albanese ($\cF^{1\ \vee}/\cF_\ZZ^\vee$) 
 varieties. The latter is related to $d\Phi$, as follows.
The infinitesimal deformations of $\cH_o$ are controlled by the Kodaira--Spencer map
$\kappa\in  T_{B,o}^\vee\otimes H^1(\cH_o,T)$.
On the other hand,
 $T_{\bD,\widetilde{\Phi}(o)}\subset H^0(\cH_o,\Omega^1)^\vee\otimes H^1(\cH_o,\cO)$,
 so cup product and   
 contraction $T\otimes \Omega^1\to \cO$ give a map $m^\vee:H^1(\cH_o,T)\to T_{\bD,\Phi(o)}$,
$m^\vee(\xi)(\alpha):=\xi\cup\alpha\in H^1(\cH_o,\cO)$.
\begin{thm}[\cite{griffiths_periods}]\label{griffiths_derivative}
The infinitesimal period map satisfies
\[
 \xymatrix{T_{\cB,o}\ar[dr]_-{\kappa}\ar[rr]^-{d\widetilde{\Phi}_o}& & T_{\bD,\widetilde{\Phi}(o)}\subset H^0(\cH_o,\Omega^1)^\vee\otimes H^1(\cH_o,\cO)\\
	& H^1(\cH_o,T)\ar[ur]^-{m^\vee} & \\},
\]
that is, 
 $\forall Y\in T_{\cB,o}$, one has
\[
 d\widetilde{\Phi}_o(Y)= \textrm{gr }\nabla^{GM}_Y=\kappa(Y)\cup\quad \in T_{\bD,\widetilde{\Phi}(o)}\subset\sym^2 H^0(\Omega^1)^\vee.
\]
\end{thm}
	\subsection{Lagrangian structures} 
 The  spaces
$\Higgs_{G,D,c}$
  are Lagrangian fibrations, and
the existence of such
 structure on a family of complex tori is not automatic,
  as shown in \cite{donagi_markman_cubic}.
Indeed, if  $h:\cH\to (\cB,o)$  carries a holomorphic symplectic form and has
Lagrangian fibres (\emph{Lagrangian structure}), then the symplectic form  induces an isomorphism $i: \cV^\vee\to T_{\cB}$.
 For a family of abelian varieties, a choice of marking gives a trivialisation  
  $\cV\simeq V\ctimes  \cO_\cB$,   $V:=H^0(\cH_o,\Omega^1)^\vee$, which,
combined with $i$,  gives
an ``affine structure''
$ \alpha: V^\vee\ctimes \cO_\cB\simeq T_\cB.$
Having Lagrangian fibres requires a compatibility between the affine  structure and the period map:
 the 
\emph{global cubic condition}.
\begin{thm}[\cite{donagi_markman_cubic},\cite{donagi_markman}]
For a
 family $h:\cH\to (\cB,o)$ of marked abelian varieties,  with a period map $\Phi$ and
 an  affine structure $\alpha:V^\vee\ctimes\cO_\cB\simeq T_{\cB}$  the following are equivalent:
			\begin{enumerate}
 			 \item $h$ admits a  Lagrangian structure
			  inducing $\alpha$
			 \item $\widetilde{\Phi}: (\cB,o)\to \sym^2 V$ is locally the Hessian of a function $\bF$
			  \item the differential
\[
 d\Phi\circ \alpha\in H^0(V\otimes \sym^2 V\ctimes \cO_\cB)\simeq  \mhom (\cV^\vee, \sym^2\cV)
\]
is the image under $\sym^3 V\hookr V\otimes \sym^2 V$ of  a cubic $c\in\sym^3 V\ctimes H^0(\cB,\cO_\cB)$.
 			\end{enumerate}
In case any of these conditions holds, there exists a unique symplectic  structure, for which the 0-section is Lagrangian.
\end{thm}
The key observation  is that if, for a contractible open $\cU\subset \cB$, the canonical symplectic form on $T^\vee_\cU$ descends to
 $\left. \cH\right|_\cU\simeq \cV/\cF_\ZZ^\vee$, translations by the family of lattices $(\Delta_\delta, Z)$ should be symplectomorphisms, 
which forces $Z$ to be a Hessian.
In terms of local coordinates and a basis of $V$, the cubic $c$ is given by 
$c=\sum\frac{\partial^3 \bF}{\partial a_i\partial a_j\partial a_k}da_i\cdot da_j\cdot da_k$
and the cubic condition is  the equality of mixed partials.
 The (non-unique)  function $\bF$
is 
called \emph{  holomorphic pre-potential} and plays a key  r\^ole  in Seiberg--Witten theory and 
special K\"ahler geometry.

\medskip

\section{Cameral Covers}\label{covers}
 
We  describe now a family of ramified Galois covers
of $X$,  with Galois group  $W=N_G(T)/T$.  
We  
 refer the reader to \cite{hitchin_sb,donagi_spectral_covers,faltings,scognamillo_elem,don-gaits}
  for more details and insight. What we call a ``cameral cover'' is an $\kd$-valued cameral cover
in the sense of \cite{donagi_spectral_covers} and \cite{don-gaits}.
      \subsection{The adjoint quotient}\label{adjoint}
 The  group $W$ is  finite,  acting  on  $\ft=\textrm{Lie}(T)$ by reflections
with respect to root hyperplanes, so by
 \cite{chevalley}, Theorem A the subalgebra of invariants
 is isomorphic to a polynomial
algebra on $l=\rk \fg$ generators:
 \[
 \CC[\ft]=\textrm{Sym } \ft^\vee \supset  \CC[\ft]^W\simeq \CC[I_1,\ldots,I_l],
\]
and  so there exists a quotient morphism $\chi:\ft\to \ft/W$.
The quotient 
$\ft/W$ is a cone: a variety with a  $\CC^\times$-action, induced by the homothety $\CC^\times$-action on $\ft$,
and   a unique fixed point, the orbit  of the origin.
The morphism $\chi$ is $\CC^\times$-equivariant.
There is an isomorphism $\ft/W\simeq \CC^l$,  but it is  non-canonical, since 
there is no preferred linearisation of the action. If we \emph{fix}  generators
$\{I_i\}$ (and positive roots $\cR^+$), then we can identify $\chi$ with a map $\ft\to\CC^l$ (or respectively, $\CC^l\to\CC^l$),
 $v\mapsto (I_1(v),\ldots, I_l(v))$.
For classical Lie algebras 
$\chi$  maps  a diagonal matrix to the
coefficients of its characteristic polynomial.

 By construction,   $\chi$ is a branched
$|W|:1$ cover, ramified along the root hyperplanes.
  It is a classical result 
(\cite{steinberg})
that 
 $\det (d\chi)=c\prod_{\alpha\in\cR}\alpha$,
for some  $c\neq 0$. 
We  call  $\fD_\chi=\prod_{\alpha\in\cR}\alpha$  
\emph{the discriminant} of $\fg$.
Since 
$\fD_\chi\in \CC[\ft]^W$,   we can write (non-canonically)
$ \fD_\chi = P(I_1,\ldots,I_l)$,
for some (weighted-homogeneous) polynomial $P$. Then $\{P=0\}\subset \CC^l$
is the equation of the
branch locus
$\bra(\chi)\subset\ft/W$, i.e., the (singular) discriminant hypersurface of $\chi$.

The inclusion $\ft\subset \fg$ induces $\CC[\ft]^W\simeq \CC[\fg]^G$,
where $G$  acts on $\fg$ by the adjoint action. Consequently, $\ft/W\simeq \fg\sslash G$ (\emph{the adjoint quotient}), and
 we can choose the $W$-invariant polynomials on $\ft$ to be restrictions of $G$-invariant polynomials
on $\fg$. In particular, $\fD_\chi$ is the restriction of a $G$-invariant polynomial, whose  locus
of non-vanishing
is the set of regular elements in $\fg$.

The degrees $\deg I_i=d_i$, or sometimes their shifts $m_i=d_i-1$,   are  called
the \emph{exponents} of $\fg$, and are important  invariants. In particular,
%
$\prod_i d_i=|W|,\quad \sum_i m_i= |\cR^+|,\quad \sum_i(2m_i+1) = \dim \fg$.
We note that 
$|\cR|=2|\cR^+|= \dim \fg-\rk \fg$ and that
$d_i\geq 2$, as $G$ is semi-simple and there is no linear invariant. 
\\
  \subsection{Cameral covers} \label{cameral}
The  $\CC^\times$-equivariance of $\chi:\ft\to\ft/W$ allows us to twist it by any 
principal $\CC^\times$-bundle, obtaining a quotient
morphism between  the associated fibre bundles.
 Applied to  
$L=K_X(D)$, this gives 
\[
\xymatrix@1{ p:& \tot \ \ft\ctimes \kd\ar[r]& \tot \ft\ctimes \kd/W},
\]
or, non-canonically, 
$\tot \kd^{\oplus l} \to \tot \bigoplus_i \kd^{d_i} $,
$s=(s_1,\ldots,s_l)\mapsto (I_1(s),\ldots, I_l(s))$.
 We denote by $\cB$ the Hitchin base, i.e., the space of global sections of the quotient
(cone) bundle
\[
 \cB=H^0(X,\ft\ctimes \kd/W)\simeq \bigoplus_i H^0(X,  \kd^{d_i}).
\]
Varying the
evaluation morphisms
$\ev_b:X\to \tot \ft\ctimes \kd/W $ with $b\in \cB$,  we get a morphism 
\[
\xymatrix@1{  \cB\times X\ar[r]^-{\ev}&\  \tot \ft\ctimes \kd/W},
\]
which is the evaluation morphism of the tautological section over 
 $\cB\times X$.
We can then pull back  the  $W$-cover $p$ to $\cB\times X$:
\[
\xymatrix{\widetilde{X}_b\ar@<-0.5ex>@{^{(}->}[r]\ar[d]^-{\pi_b}&\cX\ar[rr]^-{\io}\ar[d]^-{\pi}& &\tot \ft\ctimes \kd\ar[d]_-{p}   \\
	  \{b\}\times X \ar@<-0.5ex>@{^{(}->}[r]&\cB\times X\ar[rr]^-{\ev}\ar[dr]&& \tot \ft\ctimes \kd/W\ar[dl]^-{q}  \\
	    & & X&   \\
}.
\]
The result is  a family  $f=\textrm{pr}_1\circ \pi:\cX\to \cB$ of (not necessarily smooth, integral or reduced) covers of $X$,
with smooth  total space $\cX$, which  is  generically a $|W|:1$  cover of $\cB\times X$.
We call the fibres
 $f^{-1}(b)=\widetilde{X}_b\subset \tot \ft\ctimes \kd $ \emph{cameral covers} of $X$, and  $\pi: \cX\to \cB\times X$
\emph{the universal cameral cover}.  

We can show  that $\cX=Z(s)$, for $s$  a section of 
the pullback of $ \ft\ctimes \kd/W$
 to  the triple product $M= \cB\times X\times \tot \ft\ctimes \kd $. 
Indeed, let $r=p\circ q: \ft\ctimes \kd\to X$ be the bundle projection. The   tautological section 
$\tau$ of $q^\ast\left(\ft\ctimes \kd/W\right)$ can be pulled to
$\tot \ft\ctimes \kd$, giving a section $p^\ast \tau$ of $r^\ast\ \ft\ctimes \kd/W$.
Then $\cX=Z(s)\subset M$,
for
$ s=  \textrm{pr}_{12}^\ast\ \ev - \textrm{pr}_{3}^\ast\ p^\ast \tau \in H^0\left(M, r^\ast\ \ft\ctimes \kd/W\right)$.
        \subsection{Genericity and Discriminants}\label{discriminants}
If $\kd^2$ is very ample,
by Bertini's theorem (\cite{scognamillo_elem}, Section 1), there exists a
 Zariski open subset $\scB\subset \cB$ of  cameral covers  which are \emph{generic}, i.e., 
$\textrm{a)}\ $ are smooth and $\textrm{b)}\ $ have simple ramification.
A more concrete description of $\scB$ can be given in terms of discriminant loci.
A root $\alpha\in\ft^\vee$ induces a bundle map $r^\ast \ft\ctimes\kd\to r^\ast \kd$, so after
pullback by
the tautological section of $r^\ast \ft\ctimes \kd$,
 can be identified with a section
$\alpha\circ\io\in H^0(\cX,r^\ast \kd)$. Its  restriction  to $\widetilde{X}_b$
is  $\alpha\circ \io_b\in H^0(\widetilde{X}_b,\pi_b^\ast \kd)$.
 Consequently, $\fD_\chi$ gives a $W$-invariant section
$\widetilde{\fD}\in H^0(r^\ast \kd^{|\cR|})^W$, and 
 $2 \ram(\pi)= Z\left( \widetilde{\fD}\right)$.  
This section   descends to a section $\fD: \ft\ctimes \kd/W\to q^\ast \kd^{|\cR|}$,
whose zero locus  $Z(\fD)$ is a singular hypersurface:
\[
\xymatrix{  \cX\ar[r]^-{\io}\ar[d]^-{\pi}& \tot \ft\ctimes \kd\ar[d]_-{p}\ar[dr]^-{\widetilde{\fD}}& \\
	   \cB\times X\ar[r]^-{\ev}& \tot \ft\ctimes \kd/W\ar[r]^-{\fD}&\tot q^\ast \kd^{|\cR|} \\
	     }.
\]
We shall occasionally
write $\fD$ and $\widetilde{\fD}$ for $\ev^\ast \fD$ and $\io^\ast\widetilde{\fD}$.
A local calculation shows that  the only possible singularities of $\widetilde{X}_b$
occur at ramification points, $\textrm{a)} \Rightarrow \textrm{b)}$, and
 ``genericity''  is equivalent to
 $\ev_b(X)\cap Z(\fD)^{sing}=\varnothing$ and
$\ev_b(X)\pitchfork Z(\fD)^{sm}$.
We have that
$  \bra(\pi)= Z\left(\ev^\ast \fD\right) \subset \cB\times X$,
and away from   $0\in\cB$, 
 $\textrm{pr}_1: \cB\times X\supset \bra(\pi)\to \cB$ is a 
ramified cover, of degree $|\cR|\deg\kd$, with $\textrm{pr}_1^{-1}(b)\simeq\bra(\pi_b)\subset X$. \emph{Its} ramification locus
is $\ev^{-1}\left( Z(\fD)^{sing}\right)$, and we have
\[
\scB\varsubsetneq \cB\backslash \left\{ \textrm{pr}_1\left( \ev^{-1}\left( Z(\fD)^{sing}\right)\right) \right\}\varsubsetneq \cB .
\]
The locus $\Delta=\cB\backslash \scB$ is often termed \emph{the discriminant locus}. For later reference,
we note that over
  simply-connected opens $\cU\subset \scB$,  the connected components of
$\bra(\pi)$ are graphs of
    holomorphic  maps $c:\cU\to X$. 

\medskip
	  \section{Abelianisation}
\label{abel}
	   
The essence of the abelianisation ideology is   that one can translate
the non-abelian data of the Higgs field into abelian (spectral) data on a cover 
of $X$. While  spectral covers suffice for classical $G$,
cameral covers seems to be better suited for the case of arbitrary  structure groups.
We discuss briefly a weak version of the abelianisation theorem, and refer  to
 \cite{hitchin_sb}, \cite{bnr}, \cite{faltings}, \cite{donagi_spectral_covers}, \cite{scognamillo_elem}, 
\cite{don-gaits}), \cite{ddp} for the details.

Let  $\pi_o :\widetilde{X}_o\to X$ be a  \emph{fixed} generic cameral cover.
By Corollary 17.8 of \cite{don-gaits},  the Hitchin fibre
$h^{-1}(o)\subset \Higgs_{G,D}$ is a torsor over a \emph{generalised Prym variety} $\Prym_{\widetilde{X}_o/X}:=H^1(X,\cT)$, for
a certain sheaf $\cT$.

We identify $\mhom(\CC^\times,T)\ztimes\Pic_{\widetilde{X}_o}$ with 
the (iso classes of) $T$-bundles on $\widetilde{X}_o$  
  via  $(\mu,\scL)\mapsto \left(\scL\backslash\{0\}\right)\times_\mu T$.
As  $W=N_G(T)/T \subset \textrm{Aut}(\widetilde{X}_o)$,
     there is 
natural $W$-action on both factors of this product.
 Next, we identify the cocharacter lattice $\cchr_G=\chr_G^\vee$  with
  $\mhom(\CC^\times, T)$, and thus 
  $\cchr \ztimes \CC^\times \simeq T$.

Every     $\check \alpha\in\crts_\fg$ determines an ideal $\chr_G(\check \alpha)=\epsilon\ZZ$, $\epsilon \geq 0$.
By inspection, $\epsilon=1$, unless  $\fg=B_l=SO_{2l+1}$ (or has it as a factor),
 and $\check\alpha$ is a long coroot. In
 that case $\epsilon=2$
and $\check\alpha/2$ is a primitive element of $\cchr_G$ (Lemma 3.3, \cite{don-pan}).

In \cite{don-gaits} the authors define two abelian sheaves $\overline{\cT}\supset \cT$ on $X$. 
  The first one is $\overline{\cT} =\pi_\ast \left(\cchr\otimes\cO^\times_{\widetilde{X}_o}\right)^W $, while  the second is
\[
X\supset U\longmapsto \cT (U):=\left\{\left.t\in \Gamma\left(\pi^{-1}(U),\cchr\otimes \cO^\times_{\widetilde{X}_o}\right)^W\ \right| \left. \widetilde{\alpha}(t)\right|_{\ram_\alpha}=1,\forall \alpha\in \cR\right\}.
\]
Here $\ram_\alpha=\ram(\pi_o)\cap\{\alpha=0\}$, and $\widetilde{\alpha}$  is the 1-psg corresponding to  $\alpha$.
Notice that $\widetilde{\alpha}(t)=\pm 1$ and $\overline{\cT}/\cT$ is a torsion sheaf. In the absence of $B_l$-factors one has $\cT=\overline{\cT}$.
The neutral components $H^1(X,\overline{\cT})^0$ and $H^1(X,\cT)^0$ are (isogenous) abelian varieties. The connected components
$\Higgs_{G,D,c}$, whenever non-empty, are torsors over $H^1(X,\cT)^0$.
	    \subsection{}
Theorem 6.4  in \cite{don-gaits} describes the above torsors.
A (generic) point $[(E,\theta)]\in h^{-1}_c(o)$ determines a morphism
$h^{-1}_c(o)\to \cchr\ztimes\Pic_{\widetilde{X}_o}$ with finite fibres (\cite{scognamillo_elem}, \cite{hurtubise_kk}).
Since $H^1(X,\overline{\cT})$ is isogenous to
\[
 \left( \cchr\ztimes\Pic_{\widetilde{X}_o}\right)^{W}\subset \cchr\ztimes\Pic_{\widetilde{X}_o},
\]
the question of
 specifying a $\Prym_{\widetilde{X}_o/X}^0$-torsor breaks, up to isogeny, into:  1) determining a coset for the
$W$-invariant $T$-bundles in $ \cchr\ztimes\Pic_{\widetilde{X}_o}$, and 2) determining a torsor
for the neutral connected component.
This is quite subtle and will not be needed later, so we refer to \cite{don-gaits} for details.

Over the generic locus $\scB=\cB\backslash \Delta$, $\Higgs_{G,D}$ is a torsor over the
relative Prym fibration $\Prym_{\cX/\scB}\to \scB$.
 As $\Higgs_{G,D}$  admits local sections over simply-connected
opens $\scU\subset \scB$, one can identify 
$\left. \Higgs_{G,D,c}\right|_{\scU}\simeq \Prym_{\cX/\scU}^0 $, whenever
$\Higgs_{G,D,c}\neq \varnothing $.

\bigskip
  
\section{Proof  of Theorems A and B}\label{proof}
	\subsection{The Symplectic Leaf}
	    
Suppose $o\in \scB\subset \cB$ is  generic (in the sense of Section \ref{cameral}), and consider the
 affine subspace
\[
 \{o\} + \cB_0 = \{o\} + \bigoplus_{i=1}^r H^0(L^{d_i}(-D))\subset \cB
\]
of codimension $(\rk \fg\deg D)$.
By its genericity, 
$\{o\}$ corresponds to a regular coadjoint orbit $\bO_{o}$, with $[o]\mapsto \overline{\bO_{o}}$  under
$\cB/\cB_0\simeq \fg^\vee_D\sslash G_D$.
Fixing  $c\in\pi_1(G)$, we get
 the closure  
$h_c^{-1}(\{o\} + \cB_0) = \overline{\cS} \subset \Higgs_{G,D,c}$ 
(if $\neq \varnothing$)
of the symplectic leaf
$\cS:=S(\bO_{o})$.
 Intersecting the $\cB_0$-coset  with the generic locus we obtain a locally closed set
\[
 \bB:= \left(\{o\}+ \cB_0\right)\cap \scB  \subset \cB
\]
and, with $h_\bB:=\left. h_c\right|_{\bB}$, a Lagrangian fibration
\[
 \xymatrix@1{  \Higgs_{G,D,c} & \supset  & \cS= h^{-1}\left(\bB \right)\ \ar[r]^--{h_\bB} &\ (\bB,o)&},
\]
 which is an
analytic family of  regularly stable Higgs bundles, having smooth cameral covers with simple Galois ramification.
The corresponding family of cameral covers $f_\bB:\left.\cX\right|_{\bB}\to (\bB,o)$  parametrises deformations of $\widetilde{X}_o$ with \emph{fixed}
 intersection
$\widetilde{X}_b\cap r^{-1}(D)\subset \tot \ft\ctimes L$.

As with any integrable system, we can associate with $h_\bB$ a  family of abelian varieties -- the family of
Albanese varieties $\textrm{Alb}(h^{-1}(b))$, $b\in\bB$, see  Section 3 in \cite{danfreed}.
By the abelianisation theorem (Section \ref{abel}), the latter is the family of 
(neutral connected components of the)
generalised
Pryms associated with the family of cameral curves. 
We then have,  for  a contractible open set $\cU\subset (\bB,o)$,  the following diagram:
 \[
  \xymatrix{\Higgs_{G,D,c}		   &    \left. \cS\right|_{\cU} \ar[dr]^-{h_\cU}\ar@<0.5ex>@{_{(}->}[l]\ar[rr]^-{\simeq}                    &           				&    \Prym^0_{\cX/\cU}\ar[dl]_-{g_\cU}                    &\cX\ar[dll]^-{f_\cU}\ar[d]_-{f}\ar@<-0.5ex>@{^{(}->}[r]& \bB\times X\times\tot\ft\ctimes L\ar[dl]_-{\textrm{pr}_1}\\
		   &					 &\cU\ar[rr] &			     &\bB             &\\	
  }.
 \]
	     The family of relative Pryms carries a holomorphic symplectic structure, which can be described in terms
of cameral curves, as we now explain. Since $\bB$ is an affine space over $\cB_0$, we have
\[
 T_{\bB,o} = \cB_0\simeq H^0\left(X, \bigoplus_i L^{d_i}(-D)\right).
\]
Let $N_i$, $i=1,2,3$, denote the normal bundles for the inclusions of $\widetilde{X}_o$ in
$\left.\cX\right|_\bB$, $\cX$ and $\tot \ft\ctimes L$, respectively. We have
\[
 T_{\bB,o} \simeq H^0(\widetilde{X}_o, N_1) \simeq H^0(\widetilde{X}_o, N_2 (-r^\ast D)) \simeq H^0(\widetilde{X}_o, N_3(-r^\ast D))^W,
\]
and all of these isomorphisms are
induced by restricting the differential of the respective projection ($f_\bB$, $f$ or $r$) to  $\widetilde{X}_o$.

On $\tot \ft\ctimes L$ there is  a canonical meromorphic 2-form  $\omega_\ft\in H^0(\ft\ctimes \Omega_{\ft\ctimes L}^2(r^\ast D))^W$, induced by the 
canonical 
symplectic form on $\tot K_X$, see e.g., \cite{hurtubise_markman_rk2}. To construct it, observe  that for a  vector space $V$,
 there is a
natural $\ft$-valued skew-form on $V\oplus \left(V^\vee\otimes \ft\right)$, given by
$\left((x,\alpha\otimes s),(y,\beta\otimes t)\right)= \alpha(y)s-\beta(x)t$. This gives rise to a $\ft$-valued 2-form
on $\tot\ft\ctimes K_X$, given locally by   $dz\wedge \left(dp\otimes \sum_i dt_i\otimes \frac{\partial}{\partial t_i}\right)$.
The
local $\cO_X$-module isomorphisms $\left. L\right|_U\simeq \left. K\right|_U$ twist this to give the meromorphic 2-form
 $\omega_\ft$.
 Alternatively, $\omega_\ft$ can be constructed 
from a Liouville form:  if $\lambda$ is the   tautological section of 
$r^\ast(\ft\ctimes L)\to\tot\ft\ctimes L$,   then $\omega_\ft=-d\lambda_o$, $\lambda_o=\left. \lambda\right|_{\widetilde{X}_o}$.  

From the tangent sequence
\[
 \xymatrix@1{0\ar[r] & T_{\widetilde{X}_o}\ar[r]& \left. T_{\ft\ctimes L} \right|_{\widetilde{X}_o}\ar[r]& N_3\ar[r]& 0}
\]
one obtains, by
contraction with $\left. \omega_\ft\right|_{\widetilde{X}_o}$,  a sheaf homomorphism $N_3\to \ft\ctimes K_{\widetilde{X}}(r^\ast D)$.
In general, this map is not an isomorphism, but induces one on invariant global sections:
\[
 H^0(\widetilde{X}_o, N_3(-r^\ast D))^W\simeq_{\omega_\ft} H^0(\widetilde{X}_o, \ft\ctimes K_{\widetilde{X}_o})^W,
\] 
see
 \cite{kjiri}, \cite{hurtubise_kk}. In fact, this is a special case of a
much more general phenomenon: as shown in \cite{kjiri}, the generalised Hitchin system satisfies the rank 2 condition of
Hurtubise and Markman, and this isomorphism  becomes a special case of Proposition 2.11 in \cite{hurtubise_markman_rk2}.

The tangent sequence of $g_\bB: \Prym_{\cX/\bB}^0\to \bB$ exhibits  $T_{\Prym^0_{\cX/\bB}}$ as an extension of
$g_\bB^\ast T_\bB$ by $T_{\Prym^0/\bB}$, which, upon restriction to $P_o=\Prym_{\widetilde{X}_o/X}^0$,
gives an extension of $T_{\bB,o}\otimes \cO_{P_o}$ by $T_{P_o}=H^1(\widetilde{X}_o,\ft\ctimes\cO)^{W}\otimes \cO_{P_o}$.
As  $T_{\bB,o}\simeq H^0(\widetilde{X}_o, \ft\ctimes K)^W$,
 at any  $\scL\in P_o$ the tangent space  $T_{\Prym^0, \scL}$ is an extension of 
a pair of Serre-dual spaces,
$H^0(\widetilde{X}_o, \ft\ctimes K )^W$ and
 $H^1(\widetilde{X}_o,\ft\ctimes\cO)^W$. 
This extension can be split,  and   the canonical
 symplectic form transported
to $T_{\Prym^0, \scL}$.  Theorem 4.1 of \cite{kjiri} implies  that under the identification
$\left. \cS\right|_{\cU}\simeq \Prym^0_{\cX/\cU}$ the symplectic form obtained by moment-map
reduction coincides with the canonical one, independent of the splitting. See \cite{hurtubise_kk}, Theorem 1.10 for the case $X=\PP^1$.
	\subsection{The Donagi-Markman cubic} We now give the first  version of our result.
\begin{proposition}\label{first_version}
 Let  $\pi_o:\widetilde{X}_o\to X$ be the cameral cover, corresponding to a generic point $o\in \bB$.
Then the Donagi--Markman cubic at $o$ for $h_\bB: \left. \cS\right|_{\bB}\to \bB$
is given by
\[
 c_o: H^0(\widetilde{X}_o, \ft\ctimes K_{\widetilde{X}_o})^W \longrightarrow \sym^2 \left(H^0(\widetilde{X}_o, \ft\ctimes K_{\widetilde{X}_o})^W \right)^\vee,
\]
\[
 c_o(\xi)(\eta,\zeta)=\frac{1}{2\pi i}\int_{\widetilde{X}_o}\kappa(Y_\xi)\cup \eta\cup \zeta,
\]
where $\kappa: T_{B,o} \to H^1(\widetilde{X}_o,T_{\widetilde{X}_o})$
is the Kodaira--Spencer map of the family $\left. f_\cU:\cX\right|_{\cU}\to \cU$, and $Y_\xi$ is the preimage of $\xi$ under the isomorphism
$T_{B,o}\simeq_{\omega_\ft} H^0(\widetilde{X}_o,\ft\ctimes K_{\widetilde{X}_o})^W $.
\end{proposition}
      \subsubsection{Proof:} Either of the maps $h_\cU$ or $g_\cU$ determines the classifying map
$\Phi: \cU\to \HH_{\dim\bB}/\Gamma$. The first-order deformations of $P_o=\Prym_{\widetilde{X}_o/X}^0$
are controlled by $H^1(P_o, T_{P_o})= H^1(P_o,\cO_{P_o})\otimes  H^1(\widetilde{X}_o,\ft\ctimes\cO)^{W}$.
The polarisation on the Prym is determined by the one on $\widetilde{X}_o$, so 
$H^1(P_o, T_{P_o})\simeq H^1(\widetilde{X}_o,\ft\ctimes\cO)^{W\otimes 2}$ and
 polarisation-preserving deformations
are contained in $\sym^2 H^1(\widetilde{X}_o,\ft\ctimes\cO)^{W}\simeq \sym^2 H^0(\widetilde{X}_o, \ft\ctimes K_{\widetilde{X}_o})^{W \vee}$,
where  the last identification uses Serre duality \emph{on $\widetilde{X}_o$}. Hence the derivative of the period map for $h_\cU$
can be identified with the derivative of the period map for the family of cameral curves $f_\bB$, and we can  apply Theorem
\ref{griffiths_derivative} to the latter. Finally, we  notice  now  the map $m^\vee$ from Griffiths' theorem is dual to the multiplication map
$m: H^0(\widetilde{X}_o, \ft\ctimes K_{\widetilde{X}_o})^{W \otimes 2}\to H^0(\widetilde{X}_o,\ft\otimes \ft \otimes K^2)\to_{\textrm{tr}} H^0(\widetilde{X}_0,K^2) $.
This is quite clear, as  one defines Serre duality as a $H^1(\widetilde{X}_0,K_{\widetilde{X}_0})$-values pairing,
followed by integration (trace) map to $\CC$. For a proof, see \cite{voisin1}, Lemma 10.22 or \cite{ACGH2},
Chapter XI, section 8.  \qed
	   We make a  linear-algebraic comment on the double-dualisation used here. 
If $V$ is a   finite-dimensional vector space, the  natural isomorphism $\mhom(V^\vee,V)=\mhom(V^{\vee\otimes 2},\CC)$
is given by 
\[
f\longmapsto \left(\alpha\otimes \beta\mapsto \beta(f(\alpha))\right).
\]
 Consequently, 
 $\mhom(V^\vee,\mhom(V^\vee,V))=\mhom(V^{\vee\otimes 3},\CC)$ is given by
\[
F \longmapsto \left(Y\otimes \alpha\otimes\beta\mapsto \beta(F(Y)(\alpha))\right).\]
We apply these to $V^\vee= H^0(\widetilde{X}_o, \ft\ctimes K_{\widetilde{X}_o})^W$.
	   Our genericity assumptions   imply that $\bra(\pi_b)\cap \textrm{Supp }D=\varnothing$
for generic $b\in \bB$.
Indeed,   cameral curves from the symplectic leaf have fixed intersection with $r^\ast D$. On the other hand
(Section \ref{discriminants}), 
  $\bra(\pi)\subset \cB\times X$ is an unbranched
covering of $\scB\subset \cB$,  with local sections  by graphs of  holomorphic maps $\scB\to X$.  
Keeping some of the  ramification points fixed corresponds to holding  some of these local sections constant, i.e.,
to a fixed part of the linear system of ramification divisors, and these  form a (proper) closed subset
of $\{o\}+\cB_0$. 
Points from $\ram(\pi_b)\cap r^\ast D$,
even if non-singular, 
 do not contribute to the
Kodaira--Spencer class: in Section \ref{proof_proper} we are going to see that   $\kappa$ is determined by  the  \emph{derivatives} of the positions of
the branch points. Hence from now on we  assume that $\bra(\pi_b)\cap \textrm{Supp }D=\varnothing$.
	\subsection{Kodaira--Spencer calculation} \label{kodaira_spencer}
We turn now to  calculating the Kodaira--Spencer map 
of the   family $f_\cU:\cX_\cU\to \cU$, $o\in\cU\subset \bB$. Recall that the (global)
Kodaira--Spencer map is the  connecting homomorphism $\varkappa:\Gamma(T_\cU)\to R^1f_\ast T_{\cX/\cU}$,
and
$\kappa = \left. \varkappa\right|_o :T_{\cU,o}\to H^1(\widetilde{X}_o,T_{\widetilde{X}_o})$.
We have the following cocycle description of $\varkappa$ from \cite{kodaira_defo}, Chapter 4.
Let $\{\scU_\delta\}$ be an (acyclic) open cover of $\cX_\cU$, $\widetilde{X}_o=\bigcup_\delta \widetilde{X}_o\cap \scU_\delta$,
on which the tangent sequence of $f_\cU$ splits.
 Fix such splittings, i.e.,
  trivialisations $\Phi_\delta: \scU_\delta\simeq \Phi_\delta(\scU_\delta)\subset\cU\times  \CC$,
with   coordinates $\{\underline{\beta}=(\beta_k)_k, z_\delta\}$, where the vertical coordinates are related by 
$z_\delta=\fii_{\delta\sigma}(\underline{\beta},z_\sigma)$. Let $Y=\frac{\partial}{\partial \beta} \in\Gamma(T_\cU)$ be a vector
field. Then a \v Cech representative for
$\varkappa(Y)$ is given by the differences of the lifts of $Y$, i.e., 
by the 1-cocycle
\[
\varkappa(Y)_{\delta\sigma }=\left.\frac{\partial \fii_{\delta\sigma}}{\partial \beta}\right|_{z_\sigma=\fii_{\sigma\delta}(\beta,z_\delta)}\frac{\partial}{\partial z_\delta}.
\]
Here  $\frac{\partial \fii_{\delta\sigma}}{\partial \beta}=\cL_Y\fii_{\delta\sigma}=Y(\fii_{\delta\sigma})$ denotes  the Lie derivative.

What we need next is   a convenient open cover of $\cX_\bB$. As we are considering deformations of
a map $\pi_o:\widetilde{X}_o\to X$ with a fixed target, this cover should be adapted to the dynamics of the
branch points.
Let $\bra(\pi_o)=\{p_1,\ldots,p_N\}$,
 with  $N=|\cR|\deg \kd$.
We consider again a contractibe open  $\cU\subset (\bB,o)$, and describe a
2-set  open cover $\bU\cup\bV=\cX_\cU$.
 The first  open  is simply $\bU:= \cX_\cU\backslash \ram(\pi)$, the complement of the ramification locus.
The set $\bV$ is a neighbourhood of $\ram(\pi)$, constructed as follows.
We fix an atlas  $\fU=\left\{(U_j,z_j) \right\}$ of $X$, where 
$U_0=X\backslash\{p_1,\ldots,p_N\}$, and $\{U_j\}$ are small non-intersecting open disks, centered
at $\{p_j\}$, and by assumption,  $\textrm{supp}(D)\subset U_0$.
Being unbranched,
$\left. \cU\times X\supset \bra(\pi)\right|_\cU\to \cU$ admits holomorphic sections over $\cU$,
given by $c_j:\cU\to X$, $j=1\ldots N$.  We assume that
$c_j(o)=p_j$, and  $c_j(\cU)\subset U_j$,
i.e., the branch points of
$\pi_b$ are contained in the open disks $U_j$ for all $b\in \cU$.
This is  always possible by the open mapping theorem and by the
simply-connectedness of $\cU$.
We thus  obtain  a (trivial) disk bundle  $\cU\times X\supset  \coprod_{j\neq 0}\textrm{graph }c_j \to \cU$
and take as our second open set
$\bV:= \pi^{-1}\left( \coprod_{j\neq 0}\textrm{graph }c_j \right)\subset \left.\cX\right|_\cU$. 
The set  $\bV$  has
$\frac{1}{2}|W||\cR|\deg L$ connected components. Recall that, for simple $\fg$, the Weyl group acts transitively on
the set of roots of fixed length. Thus,  
$\ram(\pi_o)_{p_j}=\coprod_{\alpha\in\cR_j^+}\ram_{\alpha}$, and by genericity $\left|\ram_{\alpha }\right|=|W|/2$.
Here $\cR^+_j\in\{ \cR^+, \cR^+_{shrt}, \cR^+_{long}\}$, depending on whether
$\fg$ is simply laced or not, and on the point $p_j$, if it is not.
 If $\fg$ is semi-simple and not simple, there
can be different types of fibres.  Correspondingly,
$\bV=\coprod \bV_{j\alpha}$,
$\bV_{j\alpha}=\pi^{-1}\left(\textrm{graph }c_j\right)\cap \{Z(\alpha)\}=\coprod_k \bV_{j\alpha}^k$.
The connected components $\bV_{j\alpha}^k$, $1\leq k\leq |W|/2$, are tubular neighbourhoods of the connected components of
 $\ram(\pi)$. Finally, we  have $\bU\cap \bV\subset \pi^{-1}\left(\cU\times\left(\coprod_{j\neq 0} U_j\right)\right)$.
This is indeed a good cover, since $\widetilde{X}_o\cap \bU$ is affine, and $\widetilde{X}_o\cap \bV$ is Stein.
	    
    The chosen root data and atlas  allow us to               
  split the tangent sequence of $f_\cU$ on $\bU$ and $\bV$.
To be explicit,
a choice of basis of $\cB_0\simeq H^0\left(X, \bigoplus_i L^{d_i}(-D)\right)$  identifies $\cU\subset \bB$ with nested  open subsets of $\CC^{\dim\cB_0}$, and we 
denote by $\underline{\beta}= (\beta_k)_k$
the corresponding coordinates.
Then $\{\underline{\beta}, z_j\}$ are (\'etale!) coordinates 
  on $\bU\cap\pi^{-1}\left(\cU\times U_j\right)$. 
The above choices also give (candidates for) \'etale coordinates on the
$\bV_{j\alpha}$.
Indeed, each root $\alpha\in\ft^\vee$ gives $\cX_\cU\subset\cU\times\tot(\ft\ctimes L)\to \cU\times \tot L$.
Since $\textrm{supp }D\cap U_j=\varnothing$, $\left. K(D)\right|_{U_j} = \left. K\right|_{U_j}$. The 
local coordinate
$z_j$ induces a trivialisation of $K$ and we denote by
 $\fii_{z_j}$ the composition $\left. K(D)\right|_{U_j} = \left. K\right|_{U_j} \simeq U_j\times\CC\to \CC$.
Then
\[
\xymatrix@1{\alpha_j=\fii_{z_j}\circ\alpha: & \tot\left. \ft\ctimes L\right|_{U_j}\ar[r]^-{\alpha}& \left. L\right|_{U_j}\ar[r]^-{\fii_{z_j}}&\quad \CC}
 \]
and
$\{\underline{\beta},\alpha_j \}$ are \'etale  coordinates on $\bV_{j\alpha}$. Since the cover is Galois, it is enough to  choose, for each $j$,
a simple root, say $\alpha_{1j}$, and then all other $\alpha_j$ are obtained as $s\cdot\alpha_{1j}$, $s\in W$.

We make some comments and clarifications about these coordinates.
They are only \'etale coordinates, and one could  further refine the cover  to turn them into
analytic ones. As  \'etale maps induce isomorphisms on tangent bundles, this is sufficient for our purposes.
One can understand better the local picture as follows.
 The choice of simple roots determines a basis
(of fundamental coweights) in $\ft$, and by the discussion in Sections \ref{adjoint} and \ref{cameral}
the local equations of $\left.\cX\right|_{\cU\times U_j}$
are given by a system of equations 
\[
 \left|
	  \begin{array}{l}
	   I_1(\alpha_1,\ldots,\alpha_r)=b_1(\underline{\beta},z)\\
	    \ldots \\
	   I_r(\alpha_1,\ldots,\alpha_r)=b_r(\underline{\beta},z)\\
	  \end{array}
\right. ,
\]
where  $b_i$ are  holomorphic functions and $z:=z_j$. Over (a contractible subset of) $\cU\times U_0$,
we have, on a chosen  component of $\bU$, a parametrisation $\alpha_i=g_i(\underline{\beta},z)$,
and all other  components are obtained using the $W$-action, $g_i\mapsto g_i-n_{ij} g_j$.
Over the ramification locus, say, on $\bV_{j\alpha}$, $\alpha=\alpha_1$, Weierstrass preparation theorem tells us that
\[
 \alpha^2 = (z-c(\underline{\beta}))v(\underline{\beta},z)\quad \textrm{ or } z=\alpha^2 u(\underline{\beta},\alpha)+c(\underline{\beta}),
\]
and $\alpha_i=g_i(\alpha,z)$, $i\geq 2$. Here
 $u$ and $v$ are holomorphic functions, satisfying
$v(\underline{\beta} ,c(\underline{\beta}))\neq 0$, $u(\underline{\beta},0)\neq 0$ and $v(\underline{\beta},\alpha^2 u(\underline{\beta},\alpha) +c(\underline{\beta}))= 
u(\underline{\beta},\alpha)^{-1}$.

  Next,  the cameral curves may well
have ``horizontal tangents'', so
 $\fii_{z_j}\circ \alpha$  is  an \'etale  coordinate only on a (Zariski) open subset of
$\bV_{j\alpha}$.
But for us is sufficient that
this open set contains $\ram(\pi)$, which follows from the genericity and the Lie-algebraic fact
that on $\ker\alpha\cap Z(\fD)^{sm}$, $\left. \ker d\chi\right|_{\ker\alpha}=\CC\check\alpha$,
see \cite{steinberg}.

Finally, we observe that if indeed $p_j=q_i\in\textrm{supp }D$ for some $i$ and $j$,
$D=n_iq_i+\ldots$, then the frame
$\frac{dz_j}{z_j^{n_i}}$  gives  a canonical trivialisation
of $\left. K(D)\right|_{U_j}$, induced by the local coordinate $z_j$.
      \subsection{Proof of Theorem A}\label{proof_proper}
We substitute the Kodaira--Spencer class  $\kappa_{Y_\gamma}$ in  Proposition \ref{first_version}.
Applying the cocycle description of $\kappa$ to the cover $\{\bU,\bV\}$, we see that the only contributions  
 arise from  $z_j=\fii_{j\alpha}(\underline{\beta},\alpha_j)$, i.e., from the intersections $\bU\cap\bV_{j\alpha}$.
The discriminant is 
   $\widetilde{\fD}=\prod_{\alpha\in\cR}\alpha= \pm \prod_{\alpha\in\cR^+}\alpha^2$, and so
\[
 \frac{\partial_\beta \widetilde{\fD}}{\widetilde{\fD}}=\sum_{\alpha\in\cR^+}\pi^\ast\frac{\partial_\beta \alpha^2}{\alpha^2},
\]
where the $\alpha$-th  summand has second order pole along  $\ram(\pi)\cap\bV_{j\alpha}$, and is regular elsewhere.
 Fix $j_0\neq 0$, $\alpha_{j_0}\in\cR_{j_0}^+$,
 and denote $z_{j_0}$, $\alpha_{j_0}$ and $c_{j_0}$ by $z$, $\alpha$ and $c$, respectively.
Using the local equation for $\cX_\cU\subset \cU\times U_j\times\CC $ from the previous subsection, 
we have that the
 contribution to the discrimant ratio $\pi^\ast\frac{\partial_\beta \fD}{\fD} $ is
\[
-\pi^\ast\frac{\partial_\beta c}{z-c(\underline{\beta})}+\pi^\ast\frac{\partial_\beta v}{v},
\]
where the first term has a second order pole along $\alpha=0$, while the second term is regular there. On the other hand,
we have by the implicit function theorem,  
\[
 \varkappa_{\alpha z}(Y)= \frac{\alpha}{2}\left.\pi^\ast\frac{\partial_\beta \alpha^2}{\alpha^2}\right|_{z=\fii(\beta,\alpha)}\frac{\partial}{\partial \alpha} =
-\frac{\partial_\beta c}{2\alpha u(\beta,\alpha)}\frac{\partial}{\partial \alpha} +\ldots,
\]
where the first summand has a first-order pole at $\alpha=0$ and
$\ldots$ denotes terms which are regular at $\alpha=0$. Setting $\underline{\beta}=0$ gives the cocycle representatives of
 $\kappa_{\alpha z}(Y)$.  Finally, observe that for a meromorphic function $g$ the last formula implies
\[
 \textrm{Res}_{\alpha=0}\left(\kappa_{\alpha z}\lrcorner g(\alpha)d\alpha^2\right) = 
\frac{1}{2}\textrm{Res}^2_{\alpha=0} \left(\pi^\ast\frac{\partial_\beta \fD}{\fD} g(\alpha)d\alpha^2\right).
\]
Summing over $j=1\ldots N$, $\alpha\in\cR_j^+$ and the $\ZZ/2\ZZ$-cosets in $W$    completes the proof. \qed
	\subsection{Proof of Theorem B}
Recall that the cup product on $H^0(\widetilde{X}_o,\ft\ctimes K_{\widetilde{X}_o})$
is obtained by combining the Killing form $\textrm{tr} = \sum_{\alpha\in\cR}\alpha^2\in \sym^2 \ft^\vee$  with the cup product on 
$H^0(\widetilde{X}_o, K_{\widetilde{X}_o})$. 
 With the setup and notation from the previous subsection, we have
$\alpha dz = \alpha^2 (2u+\alpha\partial_\alpha u)d\alpha + \sum_k \alpha\left(\alpha^2\partial_{\beta_k}u +\partial_{\beta_k} c\right)d\beta_k$,
and hence
\[
 \frac{1}{2}\textrm{Res}^2_{\alpha=0}\left(-\pi^\ast_o\left. \frac{\partial_\beta c}{z-c(\underline{\beta})}\right|_{\underline{\beta}=0}  g(\alpha)d\alpha^2\right) =
\textrm{Res}^2_{\alpha=0}\left( \left. \frac{\alpha(\partial_\beta\lrcorner(-d\lambda))}{\alpha(\lambda)}\right|_{\underline{\beta}=0} g(\alpha)d\alpha^2\right),
\]
since $\lambda =\sum_\alpha \alpha dz \otimes \frac{\partial}{\partial \alpha}$.
 This is the contribution to the residue from each of the connected components of
$\ram_\alpha$. We complete the proof by arguing that 
\[
 \textrm{Res}^2_p\sum_{\alpha\in\cR^+}\left(\frac{\alpha(\xi)}{\alpha(\lambda_o)} \left(\sum_{\alpha'\in\cR}\alpha'(\eta)\alpha'(\zeta)\right)  \right) =
\textrm{Res}^2_p\sum_{\alpha\in\cR}\frac{\alpha(\xi)\alpha(\eta)\alpha(\zeta)}{\alpha(\lambda_o)}.
\]
Indeed, the terms with $\alpha\neq \alpha'$ vanish due to the $W$-invariance of  $\xi$, $\eta$, $\zeta$,  as follows. 
The $W$-action on $H^0(\widetilde{X}_o,\ft\otimes K)$ is a combination of the actions on $\ft$ and $\widetilde{X}_o$.
The choice
of simple roots gives, dually, a basis of $\ft$, so we have $H^0(\widetilde{X}_o,\ft\ctimes K)\simeq  H^0(\widetilde{X}_o, K)^{\oplus l}$.
Then 
$\eta=(\eta_1,\ldots,\eta_l)$ and invariance with respect to the symmetry $s_\alpha$   means
$s_\alpha^\ast\eta = s_\alpha^{-1}(\eta_1,\ldots,\eta_l)$. The Weyl group acts by $s_\alpha(\beta)=\beta-n\alpha$, hence locally
$\left. s_\alpha^\ast\eta_i\right|_{\alpha=0}=\eta_i$, and $\eta_i=f d\alpha$, for some \emph{odd} function of $\alpha$, which does not contribute
to the residue at $\alpha=0$. \qed

\bigskip

\bibliographystyle{alpha}
\bibliography{biblio}
%

\end{document}